\documentclass[10pt]{article}
\usepackage{latexsym}
\usepackage{amsfonts}
\usepackage{enumerate}
\usepackage{multicol}
\usepackage{graphicx}
\usepackage{amssymb}
\usepackage{amsmath}
\usepackage{epic}

\title{Additive Combinatorics in Groups and Geometric Combinatorics on Spheres}

\author{B\'{e}la Bajnok  \\
{\small Department of Mathematics, Gettysburg College, U.S.A.}  \\ {\small Email: bbajnok@gettysburg.edu.}}

\date{May 11, 2022}

\newtheorem{thm}{Theorem}
\newtheorem{defin}[thm]{Definition}

\newtheorem{prop}[thm]{Proposition}
\newtheorem{conj}[thm]{Conjecture}

\newtheorem{prob}[thm]{Problem}

\begin{document}

\maketitle

\begin{abstract} [This paper is based on the first of two plenary talks given by the author at the 53rd Southeastern International Conference on Combinatorics, Graph Theory \& Computing, held at Florida Atlantic University,  March 7--11, 2022.]  We embark on a tour that takes us through four closely related topics: the dual concepts of independence and spanning in finite abelian groups and the analogous dual concepts of designs and distance sets on spheres.  We review some of the main known results in each area, mention several open questions, and discuss some connections among these four interesting topics.

\end{abstract}

\thispagestyle{empty}

\section{Introduction}

In this paper we discuss four interrelated topics: two from additive combinatorics, and two from geometric combinatorics.  Namely, we introduce the concepts of {\em independence} and {\em spanning} in the context of finite abelian groups, and we review the seemingly unrelated objects of {\em designs} and {\em distance sets} on the sphere.  As we shall see, strong connections exist among these concepts: independence and spanning are complementary notions in groups, as are designs and distance sets in Euclidean space (or association schemes more generally), and a definite parallel can be established between the former concepts and the latter ones.  We will also explain how results from one of these topics can be profitably applied in the development of another.

Our itinerary for visiting these four topics is summarized as follows:

\begin{center}

\begin{tabular}{l | c c c |}

& maximize $t$ & & minimize $s$ \\ \hline 

&&& \\
{\bf Additive combinatorics} & $t$-independent &  & $s$-spanning \\

in finite abelian groups &  sets && sets \\ 
&&& \\

& $\boldsymbol{\downarrow}$ && $\boldsymbol{\uparrow}$ \\ 
&&& \\

{\bf Geometric combinatorics} & $t$-designs & $\boldsymbol{\rightarrow}$ & $s$-distance  \\

on the sphere  &  & & sets \\ 
&&& \\ \hline

\end{tabular}
\end{center}

We can introduce these four concepts heuristically, as follows.  
Given a finite subset $A$ in an abelian group $G$, we measure the degree to which $A$ is independent in $G$ by a nonnegative integer $t$, and
the efficiency with which $A$ generates $G$ by a nonnegative integer $s$.
Analogously, given a finite subset $X$ on the surface of a sphere $S$, we measure
the degree to which $X$ is well balanced on $S$ by a nonnegative integer $t$, and the number of different distances that $X$ possesses by a nonnegative integer $s$.  
We provide precise definitions below.

\section{The independence number of a subset of an abelian group} \label{FAU section 2}

Throughout, we let $G$ be an additively written abelian group of finite order $n$.   When $G$ is cyclic, we identify it with $\mathbb{Z}_n=\mathbb{Z}/n\mathbb{Z}$; we consider $0,1,\dots,n-1$ interchangeably as integers and as elements of $\mathbb{Z}_n$.  
  
Our goal is to introduce a notion for the degree to which a subset $A=\{a_1, \dots, a_m\} $ of $G$ is independent in $G$.  
Recall that, in an $R$-module $M$, a subset $A=\{a_1, \dots, a_m\} $  is {\em linearly independent} in $M$, if the only way to have $$\lambda_1 \cdot a_1 + \lambda_2 \cdot a_2 + \cdots + \lambda_m \cdot a_m =0$$ with $\lambda_1, \lambda_2, \ldots, \lambda_m \in R$ 
is to have $$\lambda_1 = \lambda_2 = \cdots = \lambda_m = 0.$$  
Clearly, our $G$ is a $\mathbb{Z}$-module, but finite, so no set $A$ satisfies this.  
We make the following definition.
  
\begin{defin} \label{FAU indep}

A subset $A$ is {\em $t$-independent} in $G$ for a nonnegative integer $t$, if the only way to have $$\lambda_1 \cdot a_1 + \lambda_2 \cdot a_2 + \cdots + \lambda_m \cdot a_m =0$$ with $\lambda_1, \lambda_2, \ldots, \lambda_m \in \mathbb{Z}$ 
  and  $\sum_{i=1}^m |\lambda_i| \le t$ 
is to have $$\lambda_1 = \lambda_2 = \cdots = \lambda_m = 0.$$

\end{defin}

As an example, consider $A=\{1,4,6,9,11 \}$ in $G = \mathbb{Z}_{25}$.  
We see that, for instance, $1+1+1+1-4$, $6+6-1-11$, and $1+4+9+11$ each equal $0$, but no such expression involving fewer than four terms is possible, so $A$ is 3-independent in $G$.  

It will be helpful to introduce the following notation: for a positive integer $h$, we let  
the {\em $h$-fold signed sumset} of $A=\{a_1, a_2, \ldots, a_m \}$ be 
$$h_{\pm} A= \left \{\lambda_1 \cdot a_1 + \cdots + \lambda_m \cdot a_m \; : \; \lambda_i \in \mathbb{Z}, \; \Sigma_{i=1}^m | \lambda_i | = h \right\}.$$ We can thus rephrase Definition \ref{FAU indep} to say that $A$ is $t$-independent in $G$ if $\cup_{h=1}^{t} h_{\pm} A$ does not contain $0$.
 Returning to our previous example, we find that 
\begin{eqnarray*} 
1_{\pm} A & = &  \{1, 4, 6, 9, 11, 14, 16, 19, 21, 24\} \\
2_{\pm} A & = & \{2, 3, 5, 7, 8, 10, 12, 13, 15, 17, 18, 20, 22, 23 \} \\
3_{\pm} A & = & \{1, 2, 3, 4, 6, 7, 8, 9, 11, 12, 13, 14, 16, 17, 18, 19, 21, 22, 23, 24 \}
\end{eqnarray*} 
However, $4_{\pm}A = \mathbb{Z}_{25}$, so $A$ is 3-independent but not 4-independent in $G$.  (It is worth noting that our chosen illustration has some special properties; for example, that $0_{\pm} A, 1_{\pm} A,$ and $2_{\pm} A$ form an exact partition of the group.)  

We should point out that independence, as defined above, has strong connections to some classical and well known concepts in additive combinatorics.  Namely: 
\begin{itemize}
  \item 
$A$ is a {\em zero-$h$-sum-free set} for some $h \in \mathbb{N}$ if
\begin{center} $x_1+x_2 + \cdots + x_h = 0$ \end{center} has no solutions in $A$.

\item
 $A$ is a {\em $(k,\ell)$-sum-free set} for some $k, \ell \in \mathbb{N}$, $k \neq \ell$ if
\begin{center} $x_1+x_2 + \cdots + x_k = y_1 + y_2 + \cdots + y_{\ell}$ \end{center} has no solutions  in $A$.  
A $(2,1)$-sum-free set is simply called a {\em sum-free set}.

\item
 $A$ is a {\em $B_h$ set} for some $h \in \mathbb{N}$ if
\begin{center} $ x_1+x_2 + \cdots + x_h = y_1 + y_2 + \cdots + y_h$ \end{center} has only trivial solutions  in $A$.  
A $B_2$ set is called a {\em Sidon set}.

\end{itemize}
  
Observe that $A$ is $t$-independent in $G$ if, and only if, for every $k, \ell \in \mathbb{N}_0$ with $k+ \ell \le t$,
$$x_1+x_2 + \cdots + x_k = y_1 + y_2 + \cdots + y_{\ell}$$ 
only has trivial solutions in $A$ ($k=\ell$ and terms are the same).  
Therefore, we may employ the classic terminology to restate Definition \ref{FAU indep} as follows:

\begin{prop}
A subset $A$ of $G$ is $t$-independent in $G$ if, and only if,  all of the following hold: 
\begin{itemize}
  \item 
$A$ is zero-$h$-sum-free for $1 \leq h \leq t$;

\item
 $A$ is $(k,\ell)$-sum-free for $1 \leq \ell < k \leq t- \ell$;

\item
 $A$ is a $B_h$ set for $2 \leq h \leq \lfloor t/2 \rfloor$.

\end{itemize}

\end{prop}
It is enough, in fact, to require these conditions for equations containing a total of $t$ or $t-1$
terms; therefore the total number of equations considered can be reduced to $2+(t-2)+1 = t+1$.

A main question then regarding $t$-independent sets in abelain groups is to see how large they can get.

\begin{prob}

For each abelian group $G$ and positive integer $t$, find the size 
$\tau (G,t)$ of the largest $t$-independent set in $G$.
\end{prob}

Since $A \subseteq G$ is $1$-independent whenever $A$ does not contain $0$, we have $\tau (G, 1) = |G|-1$ for every $G$; and since 
$A \subseteq G$ is $2$-independent if, and only if, $A$ is asymmetric (that is, $A$ and $-A$ are disjoint), we get 
$$\tau (G,2) = \frac{|G| - |L| }{2},$$ where $L$ is the set of involutions in $G$.  
  
We can evaluate $\tau (G, 3)$ in cyclic groups, as follows.  Note that $A$ is 3-independent in $G$ if, and only if, it is asymmetric, sumfree, and no three of its (not necessarily distinct) elements add to $0$.  Therefore, we can observe that the `odd' numbers $1, 3, 5, \ldots$ in $\mathbb{Z}_n$ are 3-independent as long as the largest odd element is less than $n/3$; in fact, we can go up to just below $n/2$ when $n$ is even.  We can do better when $n$ has a divisor $d$ that is congruent to 5 mod 6: the union of the first $(d+1)/6$ `odd' cosets of a subgroup of index $d$ form a 3-independent set as well.  Using Kneser's Theorem (cf.~\cite{Kne:1953a}, \cite{Nat:1996a}), we can prove that we cannot do better:

\begin{thm}[with Ruzsa, cf.~\cite{BajRuz:2003a}]
For any positive integer $n$ we have 
$$\tau(\mathbb{Z}_n,3)= \left\{
\begin{array}{cl}
\left\lfloor \frac{n}{4} \right\rfloor & \mbox{if $n$ is even;}\\ \\
\left(1+\frac{1}{p}\right) \frac{n}{6} & \mbox{if $n$ is odd, has prime divisors 5 mod 6, and $p$ is the smallest;} \\ \\
\left\lfloor \frac{n}{6} \right\rfloor & \mbox{otherwise.}
\end{array}\right.$$
\end{thm}
A more general result in \cite{BajRuz:2003a} evaluates $\tau(G,3)$ for all groups $G$, with the exception of those whose exponent is a product of a power of 3 and some primes that are congruent to 1 mod 6.  

For $t=4, 5$, and 6, we have the following computational data:
\begin{eqnarray*}
\tau(\mathbb{Z}_n, 4) & = & \left\{
\begin{array}{lll}
1 & \mbox{iff} & n \in [5,12]\\
2 & \mbox{iff} & n \in [13,26]\\
3 & \mbox{iff} & n \in [27,45], n=47\\
4 & \mbox{iff} & n=46, n \in [48,68], n=72, 73\\
5 & \mbox{iff} & n = 69, 70, 71, n \in [74,102];
\end{array}\right. \\ \\
\tau(\mathbb{Z}_n, 5) & = & \left\{
\begin{array}{lll}
1 & \mbox{iff} & n \in [6,17], n=19, 20\\
2 & \mbox{iff} & n=18, n \in [21,37], [39,41], [43,45], 47\\
3 & \mbox{iff} & n=38, 42, 46, n \in [48,69], [71,73], [75,77], 79, 81, 83, 85, 87;
\end{array}\right.\\ \\
\tau(\mathbb{Z}_n, 6) & = & \left\{
\begin{array}{lll}
1 & \mbox{iff} & n \in [7,24]\\
2 & \mbox{iff} & n \in [25,69]\\
3 & \mbox{iff} & n \in [70,151], [153,160].
\end{array}\right.
\end{eqnarray*}
(Here $[k, \ell]$ is the set of integers between $k$ and $\ell$, inclusive.)  We can prove the following bounds:

\begin{thm}[with Ruzsa, cf.~\cite{BajRuz:2003a}]

For all $t \geq 2$, $\epsilon > 0$, and large enough $n$, we have
$$\left( \frac{1}{t \lfloor (t+1)/2 \rfloor} - \epsilon \right) n^{\frac{1}{\lfloor t/2 \rfloor}} \le \tau (\mathbb{Z}_n, t) \le \left( \frac{1}{2} \cdot \left \lfloor \frac{t}{2} \right \rfloor ! \right) n^{\frac{1}{\lfloor t/2 \rfloor}}.$$
\end{thm}

Furthermore, we make the following conjecture.

\begin{conj}
Let $t \geq 2$.  The value of $$\lim \frac{\tau (\mathbb{Z}_n, t) }{ n^{1/ \lfloor t/2 \rfloor} }$$ exists if, and only if, $t$ is even.
\end{conj}


\section{Spherical designs} \label{FAU section 3}

Imagine that we want to scatter a finite number $N \in \mathbb{N}$ points on a $d$-dimensional sphere $S^d \subset \mathbb{R}^{d+1}$: how can we do it in the most ``uniformly balanced'' way?  The answer might be obvious for certain values of $N$ and $d$: arrange the points to form the vertices of a regular $n$-gon when $d=1$, or place them so that they form the vertices of a platonic solid when $d=2$ and $N \in \{4, 6, 8, 12, 20\}$.  But what should we do in general?

The answer to this question depends, of course, on how we measure the degree to which our pointset $X=\{{\bf x_1},\dots,{\bf x_N} \} \subset S^d$ is balanced.  For example, we may aim for:

\textbullet\  {\em best packing}, by maximizing 
$$\min_{i \not = j} ||{\bf x_i}-{\bf x_j}||;$$

\textbullet\  {\em best covering}, by  minimizing 
$$\max_{{\bf x} \in S^d}  \min_{i} ||{\bf x}-{\bf x_i}||;$$

\textbullet\  {\em lowest electrostatic energy (Fekete points)}, by  minimizing 
$$\sum_{i \not = j} \frac{1}{||{\bf x_i}-{\bf x_j}||}; \; \mbox{or}$$ 

\textbullet\  {\em highest thermodynamic entropy (Shub-Smale points)}, by maximizing 
$$\prod_{i \not = j} ||{\bf x_i}-{\bf x_j}||.$$
Other well-known criteria include maximizing the volume of the convex hull of $X$, or minimizing the average probability of error when signal ${\bf x_i} \in X$ is transmitted but (with some probability distribution) some ${\bf x} \in S^d$ is received instead, which then is decoded as the closest ${\bf x_j}\in X$.  
But the two concepts that have received the most attention in geometric combinatorics are: 

\textbullet\  {\em few-distance sets}: minimizing the number of distinct distances 
$$|\{||{\bf x_i}-{\bf x_j}|| : i \not = j \}|; \; \mbox{and}$$

\textbullet\  {\em spherical designs}: maximizing the degree to which the pointset is in momentum balance.  

In this section we review spherical designs and discuss how independent sets in finite abelian groups (cf.~Section \ref{FAU section 2}) are related; in Section \ref{FAU section 4} below we return to few-distance sets and their connections to the other topics of this paper.

To define spherical designs, we assume that $S^d$ is centered at the origin and has radius 1.  Given a finite pointset $X \subset S^d$ and a polynomial $f: S^d \rightarrow \mathbb{R}$, we define the average of $f$ over $X$ as 
$$\overline{f}_X = \frac{1}{|X|} \cdot \sum_{x \in X} f(x);$$ the average of $f$ over the entire sphere is given by
$$\overline{f}_{S^d} = \frac{1}{|S^d|} \cdot \int_{S^d} f(x) \mathrm{d}x$$  (here $|S^d|$ denotes the surface area of  $S^d$).  We can then make the following definition.

\begin{defin} [Delsarte, Goethals, and Seidel, cf.~\cite{DelGoeSei:1977a}] \label{FAU spherical def}

We say that a finite pointset $X \subset S^d$ is a {\em spherical $t$-design} on $S^d$ for some nonnegative integer $t$ if 
 $$\overline{f}_{X} = \overline{f}_{S^d}$$
holds for every polynomial $f: S^d \rightarrow \mathbb{R}$ of degree at most $t$.

\end{defin}

Clearly, any finite pointset $X \subset S^d$ is a spherical 0-design.  Spherical 1-designs and 2-designs have physical interpretations: it is easy to check that $X$ is a spherical 1-design exactly when $X$ is in mass balance (that is, its center of gravity is at the center of the sphere), and that $X$ is a spherical 2-design if, and only if, $X$ is in both mass balance and inertia balance. One can also verify that the vertices of a regular $N$-gon form a $t$-design on $S^1$ whenever    
$N \geq t+1$, and that the vertices of a regular tetrahedron, octahedron, and icosahedron from a $t$-design on $S^2$ for $t=2,$ $3$, and 5, respectively.  There are numerous other famous examples, such as the $24$ vertices of Neil Sloane's `improved snub cube' that form a $7$-design on $S^2$.  

The general problem is the following.

\begin{prob}
Find all positive integers $d, t,$ and $N$ for which there is a $t$-design of size $N$ on $S^d$.  
\end{prob}

There have been a number of different methods to construct spherical designs, including group theory, numerical analysis, and other techniques; for a sampling of these see \cite{Baj:1992a}, \cite{Baj:1998a}, \cite{Ban:1988a}, \cite{BanDam:1979a}, \cite{DelGoeSei:1977a}, \cite{God:1993a}, \cite{GoeSei:1979a}, \cite{HarSlo:1996a}, \cite{Hog:1996a}, 
\cite{DelGoeSei:1977a}, \cite{God:1993a}, \cite{GoeSei:1979a}, \cite{HarSlo:1996a}, \cite{Hog:1996a}, \cite{Rez:1992a}, \cite{Sei:1990a}, and \cite{SeyZas:1984a}.
Here we present a construction using additive combinatorics.

By Definition \ref{FAU spherical def}, to verify that a given pointset $X$ is a spherical $t$-design, one should check whether $\overline{f}_{X}$ agrees with $\overline{f}_{S^d}$ for every polynomial $f$ of degree at most $t$.  There is a well-known shortcut: it suffices to do this for non-constant homogeneous harmonic polynomials of degree at most $t$.  (Recall that a polynomial is called harmonic if it satisfies the Laplace equation.)  The reduction to homogeneous harmonic polynomials has two advantages.  First, there are a lot fewer of them: the set $\mathrm{Harm}_k(S^d)$ of homogeneous harmonic polynomials over $S^d$ of degree $k$ forms a vector space over $\mathbb{R}$ whose dimension is only \label{dim harm}
$$\dim \mathrm{Harm}_k(S^d) = {d+k \choose d} - {d+k-2 \choose d}.$$  Second, the average of non-constant harmonic polynomials over the sphere is zero.  Therefore, we have the following equivalent definition. 

\begin{prop} \label{FAU prop spherical}

$X \subset S^d$ is a spherical $t$-design if, and only if, $$\sum_{{\bf x} \in X} f({\bf x})=0$$
holds for all $f \in \mathrm{Harm}_k(S^d)$, $k = 1, 2, \ldots, t$.
\end{prop}

For small values of $k$ and $d$, it is possible to construct explicit bases for $\mathrm{Harm}_k(S^d)$.  In particular, for $d=1$ (and arbitrary $k$) and for $t = 1, 2, 3$ (and arbitrary $d$), we have:

\begin{tabular}{l l l}
$\mathrm{Harm}_k(S^1)$ & = & $\langle \{\mathrm{Re}(x_1+\mathrm{i} \cdot x_2)^k, \mathrm{Im}(x_1+\mathrm{i} \cdot x_2)^k \} \rangle$
\end{tabular}

\bigskip

\begin{tabular}{lll}
$\mathrm{Harm}_1(S^d)$ & = & $\langle \{x_i \; : \;  1 \leq i \leq d+1 \}\rangle$
\end{tabular}

\bigskip

\begin{tabular}{lll}
$\mathrm{Harm}_2(S^d)$ & = & $\langle \{x_i^2-x_{i+1}^2 \; : \; 1 \leq i \leq d \} \cup \{x_ix_j \; : \; 1 \leq i<j \leq d+1\}\rangle$
\end{tabular}

\bigskip

\begin{tabular}{lll}
$\mathrm{Harm}_3(S^d)$ & = & $\langle  \{x_i^3-3x_ix_j^2 \; : \; 1 \leq i \not =j\leq d+1 \} \cup \{x_ix_jx_k \; : \; 1 \leq i<j<k \leq d+1\} \rangle$
\end{tabular}

We now discuss how additive combinatorics -- in particular, $t$-independent sets in the cyclic group $\mathbb{Z}_n$ -- enables us to construct spherical $t$-designs explicitly.  

We first review the case of $d=1$, for which it is well known that the vertices of a regular $N$-gon form a $t$-design whenever $N \geq t+1$.  We identify $S^1$ with the set of complex numbers of norm 1, and set 
$X=\left \{z^j \right \}_{j=1}^{N}$,  where 
$$z = \mathrm{e}^{2 \pi \mathrm{i} /N} = \cos \left( 2 \pi  /N \right) + \mathrm{i}  \sin\left( 2 \pi  /N \right).$$ 
By Proposition \ref{FAU prop spherical}, and since $\mathrm{Harm}_k(S^1)=\langle \{\mathrm{Re}(z^k), \mathrm{Im}(z^k) \}\rangle$, $X$ is a $t$-design on $S^1$ if, and only if, $\sum_{j=1}^N (z^j)^k=0$ for all $k=1,2,\ldots, t$.  
We see that, when $k$ is not a multiple of $N$, then $z^k \neq 1$, and thus
$$\sum_{j=1}^{N} \left( z^j \right)^k = \sum_{j=1}^{N} \left( z^k \right)^j = z^k \cdot \frac{  \left( z^k \right)^{N} - 1}{z^k-1} = z^k \cdot  \frac{  \left( z^{N} \right)^k - 1}{z^k-1}   = 0.$$  Since none of $k=1,2, \ldots, N-1$ is a multiple of $N$ but $k=N$ is, $X$ is a spherical $t$-design if, and only if, $N \geq t+1$, as claimed.  (As we will soon see, there are no $t$-designs on $S^1$ of size $N \leq t$.)

We now attempt to generalize this construction for higher dimensions.  For simplicity, we assume that $d$ is odd, and let $m=(d+1)/2$.  (The case when $d$ is even can be reduced to this case by a simple technique, see \cite{Baj:1998a}.)

Let $a_1, a_2, \dots, a_m$ be integers, and set $A=\{a_1, a_2, \dots, a_m \}$.  For a positive integers $N$, define \begin{equation} \label{X(A,N)} X(A,N)=\left\{ \frac{1}{\sqrt{m}} \left( {\bf z}_N^j(a_1),  {\bf z}_N^j(a_2), \dots, {\bf z}_N^j(a_m) \right) \mbox{   } | \mbox{   }  j=1,2,\dots,N \right\},\end{equation}
where $${\bf z}_N^j(a_i)=\left( \cos \left (2 \pi j a_i /N \right), \sin \left (2 \pi j a_i /N \right) \right).$$

Note that $X(A,N) \subset S^d$.  We can then prove the following. 

\begin{thm}[\cite{Baj:1998a}] \label{FAU connect}

Let $t$, $d$, and $N$ be positive integers with $t \leq 3$, $d$ odd, and set $m=(d+1)/2$.  For integers $a_1, a_2, \dots, a_m$, define $X(A,N)$ as above.  If $A$ is a $t$-independent set in $\mathbb{Z}_{N}$, then $X(A,N)$ is a spherical $t$-design on $S^d$.  

\end{thm}

Specifically, from Theorem \ref{FAU connect}, we see that: 
\begin{itemize}
\item $X(A,N)$ is a spherical $1$-design on $S^d$ if, and only if, $A$ is 1-independent in $\mathbb{Z}_N$. For example, we may take $A=\{1\}$ for any $N \geq 2$.  (There is no spherical $1$-design on $S^d$ with $N=1$.)

\item $X(A,N)$ is a spherical $2$-design on $S^d$ if, and only if,  $A$ is 2-independent in $\mathbb{Z}_N$.
 For example, we may take  $A=\{1,2, \ldots, m \}$ for any $N \geq 2m+1 = d+2$.  (There is no spherical $2$-design on $S^d$ with $N \leq d+1$.)

\item $X(A,N)$ is a spherical $3$-design on $S^d$ if, and only if,  $A$ is 3-independent in $\mathbb{Z}_N$.
 For example (cf.~Section \ref{FAU section 2}), we may take:  
\begin{itemize} 
\item $A=\{1,3, \ldots, 2m-1\}$ with $N \geq 6m=3(d+1)$;
  \item $A=\{1,3, \ldots, 2m-1\}$ with $N \geq 4m=2(d+1)$ if $N$ is even;
\item $A=$ any $m$-subset of 
$$\left \{pi + 2j+1 \; : \; 0 \le i \le N/p-1, 0 \le j \le (p-5)/6 \right\}$$
with $$N \geq \frac{6p}{p+1} \cdot m = \frac{3p}{p+1} \cdot (d+1)$$ if $N$ is divisible by $p$ for some $p \equiv 5$ mod 6.
\end{itemize}

\end{itemize}
 Therefore, Theorem \ref{FAU connect} provides explicit constructions for all possible cases for $t=1$ and $t=2$.  For $t=3$, as we will soon see, there is no spherical 3-design on $S^d$ of size $N < 2(d+1)$.  Our construction yields all cases with even $N \geq 2(d+1)$, and all $N \geq \frac{5}{2} (d+1)$ when $N$ is divisible by 5.  In fact, with some additional techniques, we were able to prove the following.

\begin{thm} [\cite{Baj:2004a}] \label{FAU thm 3 designs}

There is no spherical 3-design on $S^d$ of size $N < 2(d+1)$.  Spherical 3-designs on $S^d$ exist when:
\begin{itemize}
\item $N$ is even and $N \geq 2(d+1)$, or
\item $N$ is odd and $N \geq \frac{5}{2} (d+1)$, but $(d,N)  \neq (2,9), (4,13)$.
\end{itemize}
\end{thm}

\begin{conj} \label{FAU conj 3-designs}

Theorem \ref{FAU thm 3 designs} gives all possible $(d,N)$ for which a spherical $3$-design on $S^d$ and of size $N$ exists.

\end{conj}
Towards Conjecture \ref{FAU conj 3-designs}, we have the following result. 
\begin{thm} [Boumova, Boyvalenkov, and Danev, cf.~\cite{BouBoyDan:2002a}]

No spherical 3-design exists on $S^d$ with $N$ odd and $$2(d+1) < N < (1+\sqrt[3]{2})(d+1) +0.300176.$$

\end{thm}
Therefore, the status quo of spherical 3-designs on $S^d$ can be summarized in the following table.

$$\begin{array}{lclcc} \hline
d & \mbox{min size} & \mbox{existence} & \mbox{non-existence} & \mbox{open}  \\ \hline \hline
1 & 4 & \geq 4 & -- & --  \\ \hline 
2 & 6 & 6, 8, \geq 10 & 7 & 9  \\ \hline
3 & 8 & 8, \geq 10 & 9 & --  \\ \hline
4 & 10 & 10, 12, \geq 14 & 11 & 13  \\ \hline
5 & 12 & 12, \geq 14 & 13 & --  \\ \hline
6 & 14 & 14, 16, \geq 18 & 15 & 17  \\ \hline
7 & 16 & 16, 18, \geq 20 & 17 & 19  \\ \hline
8 & 18 & 18, 20, \geq 22 & 19 & 21  \\ \hline
9 & 20 & 20, 22, \geq 24 & 21 & 23  \\ \hline
10 & 22 & 22, 24, 26, \geq 28 & 23, 25 & 27 \\ \hline
\end{array}$$

As we see, the smallest undecided question is whether there are nine points on $S^2$ forming a spherical 3-design; we suspect that the answer is `no'.


\section{Spherical few-distance sets} \label{FAU section 4} 

In this section we discuss the dual concept to spherical designs: spherical distance sets.  (The two topics can also be viewed in the more general setting of $P$-polynomial and $Q$-polynomial association schemes; cf.~\cite{Ban:1988a}.)

As before, we let $S^d$ be the unit sphere in $\mathbb{R}^{d+1}$, centered at the origin.  For a finite set $X=\{{\bf x_1}, \ldots, , {\bf x_N} \}$  of $N$ points on $S^d$, we let $A(X)$ be the set of distinct distances  
  $$A(X)= \{||{\bf x_i}-{\bf x_j}|| : {\bf x_i} \in X, {\bf x_j} \in X, i \not = j \}$$
that they determine.

We have the following definition.

\begin{defin}[Delsarte, Goethals, and Seidel, cf.~\cite{DelGoeSei:1977a}]

We say that a finite pointset $X \subset S^d$ is a {\em spherical $s$-code} or {\em spherical $s$-distance set} on $S^d$ for some nonnegative integer $s$ if $|A(X)|=s$.
\end{defin}

For example, it is easy to see that, for any positive integer $s$, the vertices of a regular $N$-gon form an $s$-distance set on $S^1$ for $N=2s+1$, and that the $N=d+2$ vertices of a regular simplex form a $1$-distance set on $S^d$ for any positive integer $d$.  There are many other famous and less well known examples: we just mention the regular octahedron and the regular icosahedron on $S^2$, whose vertices form a spherical 2-distance set and a spherical 3-distance set, respectively.

The central problem regarding distance sets is the following.

\begin{prob}
Find all positive integers $d, s,$ and $N$ for which there is an $s$-distance set of size $N$ on $S^d$.
\end{prob}

The following result establishes an upper bound for the size of spherical $s$-distance sets and a lower bound of the size of spherical $t$-designs, illustrating nicely the duality of the two concepts.  We set $${A(d,k)={d+\lfloor k/2 \rfloor \choose d}+{d+\lfloor (k-1)/2 \rfloor \choose d}}.$$

\begin{thm} [Delsarte, Goethals, and Seidel, cf.~\cite{DelGoeSei:1977a}] \label{FAU thm dual sph des code}

Suppose that  $X \subset S^d$ has size $N$, and that it is a spherical $t$-design as well as a spherical $s$-distance set on $S^d$.  In this case,
$${ A(d, t) \leq N \leq A(d, 2s)},$$ and thus $t \leq 2s.$
\end{thm}

And we have the following classification for the case when $t=2s$ holds.

\begin{thm} [Bannai and Damerell, cf.~\cite{BanDam:1979a, BanDam:1980a}]  \label{FAU Bannai Damerell}

Suppose that $X \subset S^d$ is a spherical $t$-design as well as a spherical $s$-distance set on $S^d$.  If $t=2s$, then $d=1$ or $s \in \{1,2\}$.  
\end{thm}

The cases of $d=1$ (regular polygons) and $s=1$ (regular simplexes) we have already mentioned above.  When $s=2$ and $t=4$, then $N = (d^2+5d+4)/2$, and we are aware of only two examples (besides the regular pentagon on the circle): the sets corresponding to $(d,N)=(5,27)$ or $(21,275)$.
  
The case of $s=2$ has received much attention even  when $t < 2s$ holds.  Denoting the maximum size of any spherical 2-distance set on $S^d$ by $M(2,d)$, we have 
$$(d^2+5d+4)/2 \leq M(2,d) \leq (d^2+3d+2)/2.$$  The lower bound  $N_d = (d^2+5d+4)/2$ is given by Theorem \ref{FAU thm dual sph des code}, while the upper bound $T_d = (d^2+3d+2)/2$ follows from the fact that the midpoints of the edges of a regular simplex on $S^d$ form a 2-distance set.   
Musin \cite{Mus:2009a} established the following results.

$$ \begin{array}{cccccc} \hline
d && N_d & T_d && M(2,d)  \\ \hline \hline
  1 && 5 & 3 &&    5 \\ \hline
2 && 9 & 6 && 6 \\ \hline
3 && 14 & 10 && 10 \\ \hline
   4 && 20 & 15 &&    16 \\ \hline
  5 && 27 & 21 &&    27 \\ \hline
6-20 && N_d & T_d && T_d \\ \hline
   21 && 275 & 253 &&    275 \\ \hline
  22 && 299 & 276 &&   276 \mbox{ or } 277 \\ \hline
23-38 && N_d & T_d && T_d \\ \hline
\end{array}$$
As we see, $M(2,s)$ agrees with $N_d$ when $d=1$, 5, or 21, and agrees with $T_d$ for all other $d$ up to $38$, except for $d=4$ and possibly $d=22$.


\section{The spanning number of a subset of an abelian group} \label{FAU section 5} 

We complete our tour through additive and geometric combinatorics by visiting our fourth topic: generating sets.  As before, we let $G$ be an abelian group of finite order $n$, written additively.  When $G$ is cyclic, we identify it with $\mathbb{Z}_n=\mathbb{Z}/n\mathbb{Z}$.
 
Our question here is as follows: how `fast' does a subset of $G$ generate the whole group (if it generates it at all)?  We make the following definition.

\begin{defin}

Let   $A=\{a_1, \dots, a_m\} $ be an $m$-subset of $G$.  We say that $A$ is an {\em $s$-spanning} set in $G$ for some nonnegative integer $s$ if every element of $G$ can be written as $$\lambda_1 \cdot a_1 + \lambda_2 \cdot a_2 + \cdots + \lambda_m \cdot a_m $$ with $\lambda_1, \lambda_2, \ldots, \lambda_m \in \mathbb{Z}$ 
 and  $\sum_{i=1}^m |\lambda_i| \le s$.

\end{defin}

Recalling  that, for a nonnegative integer $h$, the $h$-fold signed sumset of $A=\{a_1, a_2, \ldots, a_m \}$ is defined as
$$h_{\pm} A= \left \{\lambda_1 \cdot a_1 + \cdots + \lambda_m \cdot a_m \; : \; \lambda_i \in \mathbb{Z}, \; \Sigma_{i=1}^m | \lambda_i | = h \right\},$$ we can say that $A$ is an $s$-spanning set of $G$ when $\cup_{h=0}^s h_{\pm} A=G$.
  
Let us consider the (especially nice) example of $A = \{3,4\}$ in $G= \mathbb{Z}_{25}$.  Consider the following illustration. 
$$\begin{array}{ccccccc}
&&&   12 &&& \\

&&   5 &    8 &   11 && \\

&   23 &  1 &    4 &   7 &   10 & \\

  16 &   19 &   22 &    0 &   3 &   6 &   9 \\

&   15 &   18 &    21 &   24 &   2 & \\

&&   14 &   17 &   20 && \\

&&&   13 &&& 

\end{array}$$ 
In the center, we have $0_{\pm} A = \{0 \}$, the elements nearest to it form $1_{\pm} A = \{3, 4, 21, 22 \}$, the next layer is $2_{\pm} A = \{1, 6, 7, 8, 17, 18, 19, 24 \}$, and finally, the outermost elements form $3_{\pm} A = \{2, 5, 9, 10, 11, 12, 13, 14, 15, 16, 20, 23\}$.  As we can see,  $\cup_{h=0}^3 h_{\pm} A = G$, so $A$ is a 3-spanning set in $G$.  Furthermore, each element of $G$ gets generated exactly once -- more on this below.

Our main problem regarding spanning sets is the following.

\begin{prob}
For each abelian group $G$ and positive integer $s$, find the size 
$\phi (G,s)$ of the smallest $s$-spanning set in $G$.
\end{prob}

Since $A \subseteq G$ is $1$-spanning if, and only if, $\{0 \} \cup A \cup (-A) = G$, we have 
$$\phi (G, 1) = \frac{|G| + |L| -2}{2},$$ where $L$ is the set of involutions in $G$.  But $\phi (G,s)$ is not known in general for $s \geq 2$.   For $s=2$, we have the following data in cyclic groups.

$$ \phi (\mathbb{Z}_n,2)= \left\{
\begin{array}{ll}
1 & \mbox{iff $n=1, 2, 3, 4, {5};$}\\
2 & \mbox{iff $n=6, 7, \dots, 12, {13};$}\\
3 & \mbox{iff $n=14, 15, \dots, 21;$}\\
4 & \mbox{iff $n=22, 23, \dots, 33$, and $n=35;$}\\
5 & \mbox{iff $n=34$, $n=36, 37, \dots, 49$, and $n=51.$}\\
\end{array}\right.$$

We can obtain the following bounds.

\begin{prop}
For all $\epsilon, \delta$, and large enough $n$, we have
$$ \left( 1/\sqrt{2} - \epsilon \right) \sqrt{n} \le \phi (\mathbb{Z}_n, 2) \le \left(1 + \delta \right) \sqrt{n}.$$
\end{prop}

Analogously to how spherical $t$-designs and $s$-distance sets are dual concepts, so are $t$-independent sets and $s$-generating sets -- we explain this next.  
Let $a(m,s)$ denote the {\em Delannoy number} $$a(m,s)=\sum_{i\geq 0} {s \choose i} \cdot {m \choose i} \cdot 2^i.$$ 
(Delannoy numbers may also be defined by the recursion
$$a(m,s) = a(m-1,s)+a(m-1,s-1)+a(m,s-1),$$ together with the initial conditions of $a(m,0)=a(0,s)=1.$)
  
We have the following result, paralleling Theorem \ref{FAU thm dual sph des code} above.

\begin{thm} \label{FAU thm dual indep span}

Suppose that  $A \subset G$, and that it is a $t$-independent set as well as an $s$-spanning set in $G$; assume also that $t$ is even.  In this case,
$$a(m,t/2) \leq |G| \leq a(m,s),$$ and thus $t \leq 2s.$ 
\end{thm}

Cases when equality occurs in Theorem \ref{FAU thm dual indep span} are called {\em perfect}; we have the following classification:

\begin{itemize}
  \item 
$s=1$ (and $m$ arbitrary), in which case $|G|=a(m,1)= 2m+1$:  $A$ is perfect in $G$ if, and only if, $A$ and $-A$ partition $G \setminus \{0\}$.

\item

$m=1$ (and $s$ arbitrary), in which case  $|G|=a(1,s)=2s+1$: 
$A$ is perfect in $G$ if, and only if,  $G$ is cyclic of order $2s+1$, and 
$A=\{a\}$ with $\gcd (a,|G|)=1$.

\item

$m=2$ (and $s$ arbitrary), in which case  
$|G|=a(2,s)=2s^2+2s+1$: 
$A$ is perfect in $G$ if $G$ is cyclic of order $2s^2+2s+1$ and
 $A=c \cdot \{s,s+1\}$ with $\gcd (c,|G|)=1$
(For example $\{3,4\}$ in $\mathbb{Z}_{25}$ that we mentioned above).
\end{itemize}
We believe that this rather short list is complete:

\begin{conj} \label{FAU conj indep span}

The only instances of perfect sets in $G$ are those three just mentioned.  In particular, we must have $s=1$ or $m \in \{1,2\}$.

\end{conj}
It is worth noting the parallel between Theorem \ref{FAU Bannai Damerell} and Conjecture \ref{FAU conj indep span}.

For additional results and open questions on these and other related topics in additive combinatorics, we recommend the author's book \cite{Baj:2018a}.

\end{document}